\documentclass[a4paper, 11pt]{article}

\usepackage{amsmath}
\usepackage{amsfonts}
\usepackage{amsthm}
\usepackage{stmaryrd}
\usepackage{enumerate}

\usepackage[all,2cell,dvips]{xy}

\theoremstyle{plain}
\newtheorem{theorem}{Theorem}[section]
\newtheorem*{theorem*}{Theorem}
\newtheorem{lemma}[theorem]{Lemma}
\newtheorem*{lemma*}{Lemma}

\newtheorem*{prop*}{Proposition}
\newtheorem{cor}[theorem]{Corollary}
\newtheorem*{cor*}{Corollary}

\newtheorem*{example*}{Example}
\newtheorem{prob}{Problem}

\theoremstyle{definition}
\newtheorem*{defn}{Definition}

\newtheorem*{remark*}{Remark}

\newcommand{\N}{\ensuremath{\mathbb{N}}}

\usepackage{a4wide}
\usepackage{graphicx}
\usepackage{psfrag}
\usepackage{color}
\usepackage{wasysym}
\usepackage{subfigure}

\DeclareMathOperator{\pyr}{\mathsf{pyr}}
\DeclareMathOperator{\bipyr}{\mathsf{bipyr}}

\DeclareMathOperator{\cofacet}{cofacet}

\newtheorem*{balinski}{Balinski's Theorem}
\newtheorem*{gruenbaum}{Gr\"unbaum's Theorem}

\newcommand{\Pnjkm}{\ensuremath{P\left(n,j_1,k_1,\ldots,j_m,k_m\right)}}
\newcommand{\PnjkmF}{\ensuremath{\Delta_{n-1} * \left(\Delta_{j_1} \oplus \Delta_{k_1}\right) * \left(\Delta_{j_2} \oplus \Delta_{k_2} \right) * \cdots * \left(\Delta_{j_m} \oplus \Delta_{k_m}\right)}}
\newcommand{\Pnmp}[2]{\ensuremath{P\left({#1},{#2}\right)}}
\newcommand{\Pnm}{\Pnmp{n}{m}}
\newcommand{\PnmF}{\ensuremath{\Delta_{n-1} * \underbrace{\square * \cdots * \square}_{\text{$m$ times}}}}

\newcommand{\coloneqq}{\ensuremath{\mathrel{\mathop:}=}}

\definecolor{hilight}{rgb}{ .7, .7, 1 }

\newcommand{\DFGsupport}{Research supported by the \emph{Deutsche Forschungsgemeinschaft} within the research training group ``Methods for Discrete Structures''(GRK 1408)}

\begin{document}

\title{Linkages in Polytope Graphs}
\author{Axel Werner \qquad Ronald F.~Wotzlaw\thanks{\DFGsupport} \\[2ex]
\small TU Berlin, Institute of Mathematics \\
\small {\tt \{awerner,wotzlaw\}@math.tu-berlin.de}}
\date{October 19, 2007}

\maketitle

\vspace{-.2cm}

\begin{abstract}
  A graph is $k$-linked if any $k$ disjoint vertex-pairs can be joined
  by $k$ disjoint paths.
  We improve a lower bound on the linkedness of polytopes 
  slightly, which results in exact values for the minimal linkedness 
  of $7$-, $10$- and $13$-dimensional polytopes.
  
  We analyze in detail linkedness of polytopes on at most $(6d+7)/5$
  vertices. In that case, a sharp lower bound on minimal linkedness
  is derived, and examples meeting this lower bound are constructed.
  These examples contain a class of examples due to Gallivan.
\end{abstract}

\vspace{-.2cm}

\section{Introduction}

In the 2004 edition of the \emph{Handbook of Discrete and Computational Geometry}
the following question by Larman and Mani~\cite{LarmanMani1970} 
was stated as an open problem:
\begin{quote}
  \cite[Problem 20.2.6]{Kalai2004}
  Let $G$ be the graph of a $d$-polytope and $k = \lfloor d/2 \rfloor$.
  Is it true that for every two disjoint sequences $(v_1,\ldots,v_k)$ and
  $(w_1,\ldots,w_k)$ of vertices of $G$ there are $k$ vertex-disjoint
  paths connecting $v_i$ to $w_i$, $i=1,\ldots,k$?
\end{quote}
However, polytopes showing that this 
question has a negative answer in dimensions $8$, $10$, and $d \geq 12$ are not
hard to construct.
Even when $k$ is chosen as $\left\lfloor 2(d+4)/5 \right\rfloor$ such
paths do not necessarily exist. Indeed, such 
polytopes were already discovered in the 1970s by Gallivan
and later published in~\cite{McMullen1979} and~\cite{Gallivan1985}.

In customary graph theory language the above question can be rephrased as: Is the graph
of every $d$-polytope $\lfloor d/2 \rfloor$-linked? Gallivan's examples
show that this is not the case. We can then ask for the largest integer $k(d)$ such
that every $d$-polytope is at least $k(d)$-linked. This is the question of
determining minimal linkedness of polytope graphs. It is far from 
answered: Gallivan's examples show  $k(d) \leq \left\lfloor (2d+3)/5
\right\rfloor$, while a lower bound of $\lfloor (d+1)/3 \rfloor$ was proved
by Larman and Mani~\cite{LarmanMani1970}. We improve this lower bound slightly to
$\lfloor (d+2)/3 \rfloor$ in Section~\ref{sec:general}, 
which implies exact values for $k(d)$ in dimensions $7$, $10$, and $13$.

For the class of simplicial polytopes minimal linkedness 
can be determined. Larman and Mani~\cite{LarmanMani1970} 
have shown that every simplicial polytope is  at least $\lfloor (d+1)/2
\rfloor$-linked. The stacked polytopes show that this bound cannot
be improved.

\begin{figure}[ht]
   \subfigure[Simplicial $3$-polytopes are $2$-linked.]{
     \begin{minipage}[h]{.5\linewidth}
       \centering
       \includegraphics[width=.8\textwidth]{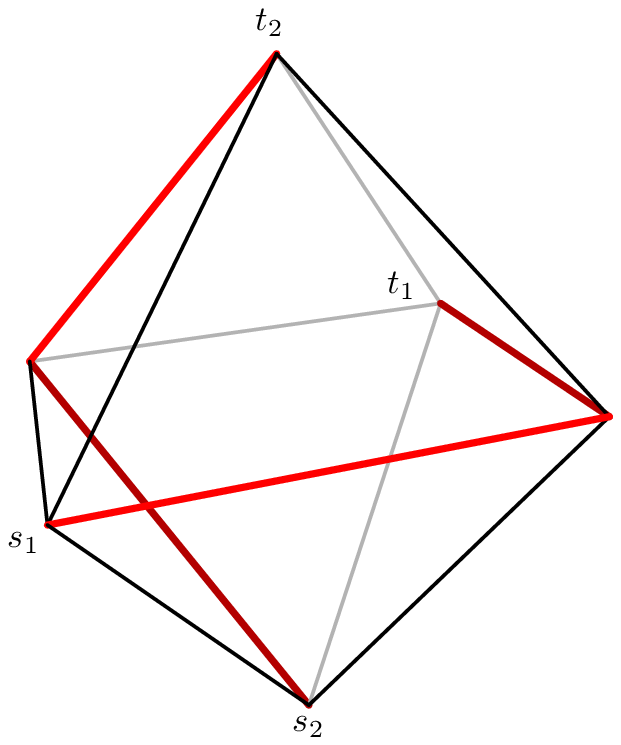}
       \label{fig:simplicial-example}
     \end{minipage}}
   \hfill
   \subfigure[Every path from $s_1$ to $t_1$ disconnects $s_2$ and $t_2$.]{
     \begin{minipage}[h]{.5\linewidth}
       \centering
       \includegraphics[width=.8\textwidth]{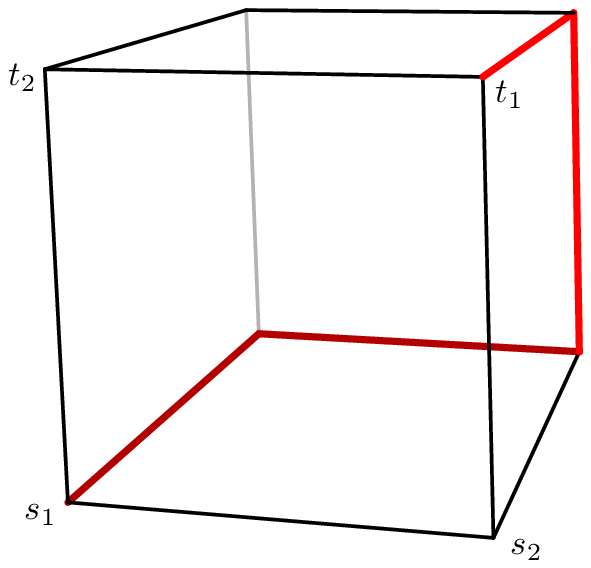}
       \label{fig:3dexample}
     \end{minipage}}
   \caption{Simplicial polytopes and $3$-dimensional polytopes.}
\end{figure}
Also, in dimensions $d \leq 5$ the values for $k(d)$ are known. While this is
trivial in dimensions $d=0,1,2$, in dimension $3$ a polytope is $2$-linked
if and only if it is simplicial and otherwise $1$-linked 
(see Figures~\ref{fig:simplicial-example} and~\ref{fig:3dexample}).
Every $4$-polytope and every $5$-polytope is $2$-linked---this follows 
from the characterization of $2$-linked graphs in~\cite{Thomassen1980} 
or the results in~\cite{Jung1970}---and examples of polytopes
that are not $3$-linked are easy to find.

Analysis of Gallivan's examples made it apparent that minimal linkedness
of $d$-polytopes on $d+\gamma+1$ vertices does depend on $\gamma$, at least if 
$\gamma$ is small. We introduce a new parameter $k(d, \gamma)$ that measures 
minimal linkedness of $d$-polytopes on $d+\gamma+1$ vertices.
We determine $k(d,\gamma)$ for polytopes on at most $(6d+7)/5$
vertices in Section~\ref{sec:few} and analyze the combinatorial types
of those polytopes with linkedness exactly $k(d, \gamma)$.
Among the combinatorial types that meet the lower bound 
Gallivan's polytopes are in some sense
the canonical ones, in some cases even unique: If $d-\gamma$
is even, there is only one combinatorial type with linkedness $k(d, \gamma)$ 
among all polytopes on $d+\gamma+1$ vertices. This type is given by
an iterated pyramid over a join of several quadrilaterals.

{\bf Acknowledgements.} The authors thank G\"unter M.\ Ziegler
for stimulating discussions on the subject, his help in preparing this paper, and
for the crucial idea in the proof of Lemma~\ref{lem:comb-char}.

\subsection*{Definitions and preliminaries}

Throughout this paper we consider polytopes only up to combinatorial equivalence,
that is, up to isomorphisms of their face lattices; none of the presented
results depend on the geometry. 
  
For $d \geq 0, \gamma \geq 0$ define the class
\[
\mathcal{P}_d^\gamma \coloneqq \{ P : P \text{ is a } d\text{-polytope on } d+\gamma+1 \text{
  vertices} \}.
\]
We denote the $d$-dimensional simplex by $\Delta_d$, the $d$-dimensional 
crosspolytope by $C_d^\Delta$ (polar to the cube $C_d$), and the
$2$-dimensional quadrilateral by $\square = C_2$.
The interested reader will find plenty of information about polytopes in
the books by Gr\"unbaum~\cite{Gruenbaum1967} and Ziegler~\cite{Ziegler1995}.

If $P$ is a polytope we denote by $\dim(P)$ the dimension of $P$ and by
$f_0(P)$ the number of vertices. We define $\gamma(P) = f_0(P) - \dim(P) - 1$.

The \emph{graph of a polytope $P$} is the graph $G(P)$
defined by the vertices and edges of $P$ and their incidence relations.

Let $G=(V,E)$ be a graph and $k \in \N$. Then $G$ is called
\emph{$k$-(vertex-)connected} if $|V| > k$ and if 
for each $C \subseteq V$ of cardinality $|C| < k$ the
graph $G - C$ is connected. The graph $G$ is \emph{$k$-linked} 
if $|V| \geq 2k$, and if for every choice of
$2k$ distinct vertices $s_1, \ldots, s_k, t_1, \ldots, t_k$
there exist $k$ disjoint paths $L_1,\ldots, L_k$ such that
$L_i$ joins $s_i$ and $t_i$ for $i=1, \ldots, k$. This implies
that $G$ is at least $(2k-1)$-connected.
The paths $L_1, \ldots, L_k$ are called a \emph{linkage} for 
$s_1, \ldots, s_k, t_1, \ldots, t_k$.

We say that a polytope $P$ is $k$-linked if the graph $G(P)$
is $k$-linked and define the following parameters:
\begin{eqnarray*}
  k(P) & \coloneqq  & \max \{ k : \text{ $P$ is $k$-linked} \} \\
  k(d, \gamma) & \coloneqq & \max \{ k : \forall P \in \mathcal{P}_d^\gamma: k(P) \geq k \} \\
  k(d) & \coloneqq & \min_\gamma k(d, \gamma)
\end{eqnarray*}

A graph $G'$ is a \emph{subdivision of $G$}
if $G'$ is obtained from $G$ by replacing each edge $uv \in E$ of $G$ by
a path $M_{uv}$ with end-vertices $u$ and $v$ (possibly of length one). 
We call the set of all interior vertices of these paths the
\emph{subdividing vertices}, the other vertices the \emph{branch vertices}.
If there is a vertex $v \in V$ such that the set of branch
vertices is $\{ v \} \cup U$ with $U \subseteq N(v)$,
we say that $G'$ is a subdivision of $G$ \emph{rooted at $v$}.

If $M$ and $N$ are paths in a graph, we write $MN$ for the union of $M$ and $N$.

Many arguments in this paper crucially depend on the following theorems
by Balinski and Gr\"unbaum.

\begin{balinski}[1961~\cite{Balinski1961}]
\label{thm:balinski}
Let $P$ be a $d$-polytope and $G=G(P)$ be its graph.
Then $G$ is $d$-connected.
\end{balinski}

\begin{gruenbaum}[1965~\cite{Gruenbaum1965}~{\cite[Section 11.1, p.~200]{Gruenbaum1967}}]
\label{thm:subdivision}
Let $P$ be a $d$-polytope, $v \in V(P)$ a vertex of $P$, and
$G = G(P)$ the graph of $P$.
Then $G$ contains a subdivision of $K_{d+1}$ rooted at $v$.
\end{gruenbaum}

The original wording of Gr\"unbaum's theorem is 
different: It is not mentioned that the subdivision can be
chosen rooted at a specified vertex. However, this extension
is an obvious by-product of Gr\"unbaum's proof.

Both theorems were proved by Barnette~\cite{Barnette1973b} for
structures more general than polytopes.

\section{Lower and upper bounds on minimal linkedness of polytopes}
\label{sec:general}

In this section we provide lower and upper bounds on $k(d)$ for
general polytopes in arbitrary dimension $d$.
We show that there is an upper bound on $k(d,\gamma)$
that is independent of $\gamma$.

\subsection{A lower bound on minimal linkedness}
\label{sec:general-lb}

Larman and Mani~\cite{LarmanMani1970} have shown that every $2k$-connected
graph that contains a $K_{3k}$ subdivision is $k$-linked. 
This statement also follows from a more
general result by Robertson and Seymour~\cite{RobertsonSeymour1995}.
Together with
Balinski's theorem and Gr\"unbaum's theorem we conclude that every
$d$-polytope is $\left\lfloor (d+1)/3 \right\rfloor$-linked.
However, already in dimension $4$ this bound is not tight. It is easy to 
see by a geometric argument and also follows from 
the characterization of $2$-linked graphs in~\cite{Thomassen1980} 
or the results in~\cite{Jung1970} that every $4$-polytope is $2$-linked.

We improve Larman and Mani's 
bound slightly by taking a closer look at the graph structure
of $d$-polytopes. 
The following argument is a variation of the proof of Larman and Mani's result
given in~\cite[pp. 70--71]{Diestel2005}.

\begin{lemma}
\label{lem:linked-subdivision}
Let $G=(V,E)$ be a $2k$-connected graph. Suppose that for every 
vertex $v$ of $G$ the graph $G$ contains a subdivision of $K_{3k-1}$ rooted at $v$.
Then $G$ is $k$-linked.
\end{lemma}

\begin{proof}
See Figure~\ref{fig:linked-subdivision} for an illustration of the proof.
\begin{figure}[ht]
  \psfrag{s1}{$s_1$}
  \psfrag{sk}{$s_k$}
  \psfrag{v1}{$v_1$}
  \psfrag{vk}{$v_k$}
  \psfrag{t1}{$t_1$}
  \psfrag{w1}{$w_1$}
  \psfrag{tk}{$t_k$}
  \psfrag{u1}{$u_1$}
  \psfrag{W'}{$W'$}
  \psfrag{W}{$W$}
  \centering
  \includegraphics[scale=.5]{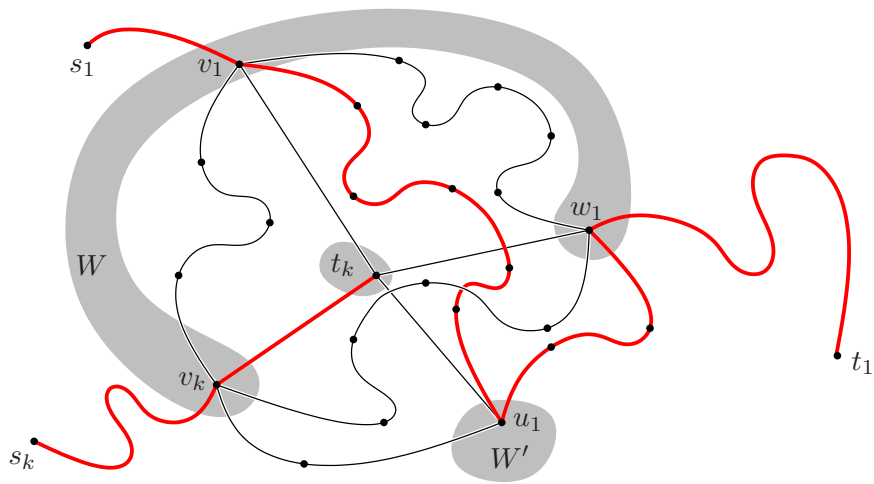}
  \caption{Illustration of the proof of Lemma~\ref{lem:linked-subdivision} with $k=2$.}
  \label{fig:linked-subdivision}
\end{figure}

Let $s_1, \ldots, s_k, t_1, \ldots, t_k$ be distinct vertices of $G$.
Let $K$ be a subdivision of $K_{3k-1}$ rooted at vertex $t_k$ with
branch vertices $ U \coloneqq \{ t_k \} \cup U'$, for $U' \subseteq N(t_k)$.

Since $G \setminus \{ t_k \}$ is $(2k-1)$-connected there exist $2k-1$ 
disjoint paths $S_1, \ldots, S_k, T_1, \ldots, T_{k-1}$ in $G$ 
avoiding $t_k$ 
such that $S_i$ joins $s_i$ to $U'$, for $i=1, \ldots, k$, and
$T_i$ joins $t_i$ to $U'$, for $i=1, \ldots, k-1$. 
Moreover, we assume that the paths have been chosen such that they
do not have interior vertices in $U'$ (and thus also not in $U$)
and that their total number of edges outside of $E(K)$ is minimal.

Let $W = \{ v_1, \ldots, v_k, w_1, \ldots, w_{k-1} \}$ 
be the vertices of these paths in $U'$, where $v_i$ is in $S_i$
and $w_i$ is in $T_i$. We then have a partition of $U$ into sets 
$\{ t_k \}$, $W$ and $W' \coloneqq U' \setminus W$ with $|W'| = k-1$.
Let $u_1, \ldots, u_{k-1}$ be the vertices in $W' \subseteq U$.
We call these vertices \emph{free}.

Since the path $S_k$ joins $s_k$ to a neighbor of $t_k$ the 
path $L_k \coloneqq S_k t_k$ joins $s_k$ and $t_k$. 

Now fix some $i \in \{ 1, \ldots, k-1 \}$ and let $M_i$ be the path in $K$
from the free vertex $u_i$ to $v_i$ and $N_i$ be the path in $K$
from $u_i$ to $w_i$.
Since the paths $S_1, \ldots, S_k, T_1, \ldots, T_{k-1}$ were
chosen minimal with respect to their number of edges outside of $K$
and $u_i$ is a free vertex, the paths $S_j$
are disjoint from $M_i$ for $j \neq i$, and they are disjoint from $N_i$ for all $j=1, \ldots, k$. 
Similarly, the paths $T_j$
are disjoint from $N_i$ for $j \neq i$, and they are disjoint from $M_i$ for all
$j=1, \ldots, k-1$. Hence we can join $v_i$ to $w_i$ via the free vertex $u_i$.

We get pairwise disjoint paths
\[
L_i = \left\{ \begin{array}{c@{\quad,\qquad}l} 
       S_i M_i N_i T_i &  1 \leq i \leq k-1\\
       S_k t_k & i=k 
    \end{array} \right. 
\]
such that $L_i$ joins $s_i$ and $t_i$, that is, a linkage for the
vertices $s_1,\ldots,s_k,t_1,\ldots,t_k$.
\end{proof}

\begin{theorem}
\label{thm:general-lower}
Every $d$-polytope is $\left\lfloor (d+2)/3 \right\rfloor$-linked.
Thus, $k(d, \gamma) \geq k(d) \geq \left\lfloor (d+2)/3 \right\rfloor$ for
all $\gamma \geq 0$ and $d \geq 1$.
\end{theorem}

\begin{proof}
Let $P$ be a $d$-polytope and $G=G(P)$ its graph.
The statement is clearly true for $d = 1$.
We set $k \coloneqq \left\lfloor (d+2)/3 \right\rfloor$.
For $d \geq 2$ we then have $d \geq 2k$ and $d+1 \geq 3k-1$.
Therefore, by Gr\"unbaum's theorem, the graph $G$ contains a
$K_{3k-1}$ subdivision at every vertex and, by
Balinski's theorem, $G$ is $2k$-connected.
By Lemma~\ref{lem:linked-subdivision}, the graph of $P$ is
$k$-linked.
\end{proof}

\subsection{An upper bound on minimal linkedness}
\label{sec:general-ub}

\begin{theorem}
\label{thm:general-ub}
Let $d \geq 2$ and $\gamma \geq 1$. The minimal linkedness of $d$-polytopes
on $d+\gamma+1$ vertices satisfies
  \[ k(d,\gamma) \; \leq \; \left\lfloor d/2 \right\rfloor. \]
\end{theorem}
\begin{proof}
For $d = 2$ the assertion is trivially true.

Let $d \geq 3$ and $\gamma \geq 1$. To prove the statement
we have to construct a $d$-polytope on $d+\gamma+1$ vertices with $k(P) \leq
\left\lfloor d/2 \right\rfloor$.

For this let $Q$ be a $3$-polytope on $4 + \gamma$ vertices that has
a square facet. For instance, for $\gamma = 1$ take the pyramid over a square
and for $\gamma > 1$ stack this pyramid $\gamma-1$ times over triangular facets.
Let $P \coloneqq \pyr^{d-3}(Q)$, the $(d-3)$-fold pyramid over $Q$. Then
$P$ is a $d$-polytope and has $d+\gamma+1$ vertices. Additionally,
$P$ is not $(\lfloor d/2 \rfloor + 1)$-linked. To see this
let $s_1, t_1, s_2, t_2$ be the vertices of a square facet of $Q$
(in that order around the facet). Then, by planarity, these cannot be linked
in $G(Q)$. Additionally, with $m = \lfloor (d-3)/2 \rfloor$ there are
$2m$ vertices left in $V(P) \setminus V(Q)$ if $d$ is odd and
$2m+1$ if $d$ is even. We choose distinct vertices $s_3,\ldots,s_{m+2}, t_3,\ldots,t_{m+2}$
arbitrarily from the set $V(P) \setminus V(Q)$ and, if $d$ is even,
$s_{m+3}$ the last vertex left in $V(P) \setminus V(Q)$ and
$t_{m+3}$ arbitrarily from $V(Q) \setminus \{ s_1, s_2, t_1, t_2 \}$.
This set of $\lfloor d/2 \rfloor + 1$ pairs of vertices cannot be linked
in $P$. Therefore $k(P) \leq \lfloor d/2 \rfloor $.
\end{proof}

In the special case $\gamma = 0$ we trivially have 
$k(d,\gamma) = \lfloor (d+1)/2 \rfloor$,
as the $d$-simplex is $\lfloor (d+1)/2 \rfloor$-linked.

Theorem~\ref{thm:general-ub} implies that $k(d) \leq \lfloor d/2 \rfloor$; 
in the next section this bound will be significantly improved.

\section{Linkages in polytopes with few vertices}
\label{sec:few}

We now study linkedness of polytopes having rather few vertices compared to their
dimension.

If we have $\gamma \leq (d+2)/5$, we can
precisely determine the value of $k(d,\gamma)$. 
However, most statements in this section make sense for all $\gamma \leq d$ but
not for $\gamma > d$. We therefore require $\gamma \leq d$ throughout the whole section.

The theory of polytopes with few vertices is closely linked to the theory
of Gale diagrams. However, we will not use Gale diagrams and prove all statements
combinatorially.

We need two basic operations that create new polytopes from given ones.
By $P_1 * P_2$ we denote the \emph{join} of two polytopes $P_1$ and $P_2$.
For example, the join of a polytope $P$ with an additional vertex $v$, 
that is, with a $0$-dimensional polytope,
results in the pyramid $\pyr P = P * v$ over $P$.
Similarly, $P_1 \oplus P_2$ denotes the \emph{sum} of the polytopes $P_1$ and $P_2$.
A special case is the sum of a polytope $P$ and an interval $I$, 
which yields the bipyramid $\bipyr P = P \oplus I$ over $P$.

\subsection{A lower bound for polytopes with few vertices}

Linkedness of a graph is a local property in the following sense:
If a graph is highly connected, then a $k$-linked subgraph ensures
$k$-linkedness for the whole graph. This is made precise in the
following lemma.

\begin{lemma}
\label{lem:subgraph}
Let $G=(V,E)$ be a $2k$-connected graph and $G'$ a subgraph of $G$ that
is $k$-linked. Then $G$ is $k$-linked.
\end{lemma}
\begin{proof}
Let $s_1, \ldots, s_k, t_1, \ldots, t_k$ be a pairing of distinct vertices
in $G$. Since $G$ is $2k$-connected, there exist $2k$ vertex disjoint 
paths $S_1, \ldots, S_k, T_1, \ldots, T_k$ such that $S_i$ connects $s_i$ to $G'$
and $T_i$ connects $t_i$ to $G'$. We choose
the paths such that each contains only one vertex from $G'$.
Let $\{ s_i' \} = G' \cap S_i$ and $\{ t_i' \} = G' \cap T_i$.
Since $G'$ is $k$-linked there exists a linkage $L_1', \ldots, L_k'$
in $G'$ for the distinct vertices $s_1', \ldots, s_k', t_1', \ldots, t_k'$ and
\[
L_i =  S_i L_i' T_i \quad, \qquad  1 \leq i \leq k
\]
is a linkage for $s_1, \ldots, s_k, t_1, \ldots, t_k$ in $G$.
\end{proof}

We obtain a lower bound on linkedness of polytopes with few vertices
by finding a highly-linked subgraph in the graph of $P$. This
highly-linked subgraph is a complete subgraph: the graph of a simplex
face of high dimension.

\begin{lemma}[\cite{Kalai1994}]
\label{lem:simplex}
Let $P$ be a $d$-polytope on $d+ \gamma +1$ vertices. 
Then $P$ has a $(d-\gamma)$-face that is a simplex.
\end{lemma}

\begin{proof}
It is easily checked that the statement is true for every $2$-polytope
on $3+\gamma$ vertices, $\gamma\geq 0$. 

Let $P$ be a $d$-polytope, $d \geq 3$. Choose a facet $F$, which is of
dimension $d'=d-1$ and has $d'+\gamma'+1$ vertices, where $0 \leq \gamma' \leq \gamma$.
By induction, $F$ has a simplex face $S$ of dimension
$\dim (S) = d'- \gamma' = d - 1 - \gamma' = d - (\gamma' + 1)$.
If $\gamma \geq \gamma' + 1$, then $\dim(S) \geq d-\gamma$ and we are done.
If $\gamma = \gamma'$, then $V(P) \setminus V(F) = \{v\}$ and $P = F * v$
is a pyramid over $F$. Hence $S * v$ is a face of $P$ and a simplex
of dimension $\dim S + 1 = d-\gamma$.
\end{proof}

\begin{theorem}
\label{thm:lb-few}
Let $d \geq \gamma \geq 0$. Then
  \[ k(d,\gamma) \; \geq \; \left\lfloor \frac{d-\gamma+1}{2} \right\rfloor.\]
\end{theorem}
\begin{proof}
  For the special cases $d=0$, $d=1$ as well as $\gamma=0$ (with arbitrary $d$)
  the assertion is trivially true.
  For $d \geq 2$, $\gamma \geq 1$ it follows directly from
  Lemmas~\ref{lem:subgraph} and~\ref{lem:simplex}, since
  $2 \lfloor (d-\gamma+1)/2 \rfloor \leq d-\gamma+1 \leq d$, and the graph of
  a $d$-polytope is at least $d$-connected, by Balinski's theorem.
\end{proof}

\subsection{An upper bound for polytopes with few vertices}
\label{subsec:few-vert-upper-bound}

To prove a good upper bound on the number $k(d,\gamma)$ we have to find a
polytope $P$ on $d+\gamma+1$ vertices with small $k(P)$.
For $\gamma \leq (d+2)/5$ 
the lower bound from Theorem~\ref{thm:lb-few} can be attained.

The class of examples we describe here was first
discovered in this context by Gallivan~\cite{Gallivan1985}, 
who constructed it using Gale diagrams.

\begin{defn}
For integers $n,m \geq 0$ and $j_1, k_1, \ldots, j_m, k_m \geq 1$ define
\[
\Pnjkm \coloneqq \PnjkmF
\]
and
\[
\Pnm \coloneqq P(n, \underbrace{1,\ldots,1}_{2m\text{ times}}) = \PnmF .
\]
\end{defn}

We consider the
complement graph of $G(\Pnm)$ to examine the linkedness of the polytopes $\Pnm$:
\begin{center}
  \psfrag{G}{\hskip -2cm $\overline{G}(\Pnm):$}
  \psfrag{ldots}{$\cdots$}
  \psfrag{egs}{$\underbrace{\text{\mbox{}\hspace{3.2cm}\mbox{}}}_{2m\text{ edges}}$}
  \psfrag{vts}{$\underbrace{\text{\mbox{}\hspace{2.1cm}\mbox{}}}_{n\text{ pyramid vertices}}$}
  \includegraphics[scale=.8]{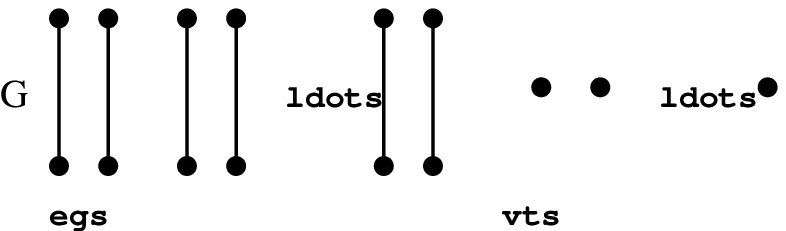}
  \bigskip
\end{center}
Roughly speaking, the reason for the low linkedness of $\Pnm$ is that
there are few vertices that can be used on a ``detour'' for a linkage
between the $m$ pairs that are not connected by an edge.

The parameters $d$ and $\gamma$ for $\Pnm$ can be determined by observing
that $\Pnm$ has $4m+n$ vertices, so $d+\gamma+1 = 4m+n$, and dimension
$\dim(\Pnm)=n-1+3m$.
Therefore we have
\begin{equation}
  \label{eqn:d-reform}
  d \; = \; n-1+3m,
\end{equation}
\begin{equation}
  \label{eqn:gamma-reform}
  \gamma \; = \; m .
\end{equation}

\begin{lemma} \label{lem:kPnm}
  Let $n,m \geq 0$ be integers. The linkedness of $\Pnm$ is given by
  \[ k(\Pnm) \; = \; \left\{ \begin{array}{c@{\quad,\qquad}l} 
      \left\lfloor \frac{4m+n}{3} \right\rfloor & n \leq 2m - 1 \\[1ex]  
      \left\lfloor \frac{2m+n}{2} \right\rfloor & n \geq 2m - 1.
    \end{array} \right. \]
If we use substitutions {\rm(\ref{eqn:d-reform})} and {\rm(\ref{eqn:gamma-reform})},
this evaluates to
  \[ k(\Pnm) \; = \; \left\{ \begin{array}{c@{\quad,\qquad}l} 
      \left\lfloor \frac{d+\gamma+1}{3} \right\rfloor & d \leq 5\gamma - 2 \\[1ex]
      \left\lfloor \frac{d-\gamma+1}{2} \right\rfloor & d \geq 5\gamma - 2.
    \end{array} \right. \]
\end{lemma}
\begin{proof}
To prove the upper bound on $k(\Pnm)$ we exhibit a ``worst possible''
pairing of the vertices of \Pnm. It is easy to see that for the
given example we have to pair as many vertices defining
an edge in the complement graph $\overline{G}(\Pnm)$ as possible. Those pairs
will necessarily block a third vertex when they are connected by a path
in $G(\Pnm)$.

If $n \leq 2m-1$, we choose $\lfloor (4m+n)/3 +1 \rfloor$ edges of 
$\overline{G}(\Pnm)$ as pairs. This many edges exist in $\overline{G}(\Pnm)$, and
to connect one of these pairs in $G(\Pnm)$ we have to use one additional
vertex. However, there are only $4m+n$ vertices altogether,
which is not enough to connect all $\lfloor (4m+n)/3 +1 \rfloor$ pairs.
This shows that $k(\Pnm) \leq \lfloor (4m+n)/3 \rfloor$ if $n \leq 2m-1$.

To show the reverse inequality we have to find a linkage for $\lfloor (4m+n)/3 \rfloor$ pairs
of vertices. Note that every pair can be connected by a path using at most one
other vertex. Also, 
each of the $4m+n-2\lfloor (4m+n)/3 \rfloor$ vertices not in the $\lfloor (4m+n)/3 \rfloor$ pairs can
be used as such a ``detour vertex.''
As we have the inequalities
\begin{eqnarray*}
  4m+n - 2 \left\lfloor \frac{4m+n}{3} \right\rfloor 
  & \geq & 3 \left\lfloor \frac{4m+n}{3} \right\rfloor - 2 \left\lfloor \frac{4m+n}{3} \right\rfloor
  \; = \; \left\lfloor \frac{4m+n}{3} \right\rfloor, 
\end{eqnarray*}
these are enough to connect all pairs in $G(\Pnm)$ by disjoint paths.

If $n \geq 2m$, we choose as pairs all $2m$ edges
in $\overline{G}(\Pnm)$ and additionally as many of the remaining
$n-2m$ isolated vertices as possible. This leaves us with
\[ 2m + \left\lfloor \frac{n-2m}{2} \right\rfloor \; = \; \left\lfloor \frac{2m+n}{2} \right\rfloor \]
pairs that can be linked with at most one more vertex remaining.

In the case $n=2m-1$ the construction described earlier applies, but still
both formulas provide the same value for $k(\Pnm)$:
\[
  \left\lfloor \frac{4m+n}{3} \right\rfloor \; = \;
  \left\lfloor \frac{6m-1}{3} \right\rfloor \; = \;
  2m - 1 \; = \;
  \left\lfloor \frac{4m-1}{2} \right\rfloor \; = \;
  \left\lfloor \frac{2m+n}{2} \right\rfloor .
\]

The lower bound for $n \geq 2m-1$ follows from Theorem~\ref{thm:lb-few}
and Equations (\ref{eqn:d-reform}) and (\ref{eqn:gamma-reform}).
\end{proof}

\begin{example*}
Let $d=8$ and $\gamma = 2$. Then $n=3$ and $m=2$ and we obtain the
$8$-polytope
\[
P \; \coloneqq \; P(3,2) \; = \; \Delta_2 * \square * \square \; = \; \pyr^3 (\square * \square).
\]
The complement of the graph of $P$ consists of $4$ disjoint edges and
$3$ isolated vertices. Obviously, $G(P)$ is not $4$-linked.
\end{example*}

In combination with Theorem~\ref{thm:lb-few} we obtain the following result.
\begin{theorem}
  Let $d \geq 0$ and $(d+2)/5 \geq \gamma \geq 0$. Then:
  \[
  k(d,\gamma) = \left\lfloor \frac{d-\gamma+1}{2} \right\rfloor.
  \]
\end{theorem}

Choosing $\gamma = \lfloor (d+2)/5 \rfloor$, we obtain Gallivan's examples, and 
the bound of the last theorem implies the following bound on $k(d)$ first
given in~\cite{Gallivan1985}.
\begin{cor}
\label{cor:gallivan-bound}
For minimal linkedness of $d$-polytopes we have
\[
k(d) \leq \left\lfloor (2d+3)/5 \right\rfloor.
\]
\end{cor}

\subsection{Analysis of polytopes meeting the lower bound}
\label{sec:char}

\begin{lemma}
\label{lem:comb-char}
Let $P$ be a $d$-polytope with the following property:
Every  facet $F$ of $P$ satisfies $|V(P) \setminus V(F)| \leq 2$.
Then $P$ is of the form
\[
\Pnjkm = \PnjkmF
\]
where $k_1 \ldots, k_m, j_1, \ldots, j_m \geq 1$ and
$d = n - 1 + j_1 + k_1 + \ldots + j_m + k_m  + m$.
\end{lemma}

\begin{proof}
The property $|V(P) \setminus V(F)| \leq 2$ implies that the
hypergraph of facet-complements, that is, the hypergraph
\[
G_{\cofacet}(P) \coloneqq (V(P), \{ W \subseteq V(P) : V(P) \setminus W \text{ is vertex set of a facet of $P$} \})
\]
is a graph (with no parallel edges, but possibly with loops).
The edges of $G_{\cofacet}(P)$ are in bijection with the facets of $P$.
Since the combinatorial type of a polytope is determined by the
vertex-facet incidences, the combinatorial type of $G_{\cofacet}(P)$ determines
the combinatorial type of $P$.

For $Q = P(n, j_1, k_1, \ldots, j_m, k_m)$ the graph $G_{\cofacet}(Q)$ 
is a disjoint union of $n$ copies of the graph that consists of one single vertex
and one single loop, and complete bipartite graphs
$K_{j_1,k_1}, \ldots, K_{j_m, k_m}$.
Thus, we have to show that $G_{\cofacet}(P)$ is of this type.
It is easy to see that loops can only occur at isolated vertices,
and that there are no vertices of degree $1$ in $G_{\cofacet}(P)$
(we follow the convention that loops contribute two edges to the degree count).
Then it suffices to check the following two properties of
$G_{\cofacet}(P)$:
\begin{enumerate}[(i)]
\item 
\label{item:prop-1}
The graph $G_{\cofacet}(P)$ does not have odd cycles. 
\item 
\label{item:prop-2}
Whenever there is a path $v_1v_2v_3v_4$ of length $3$ 
in $G_{\cofacet}(P)$, then $\{ v_1, v_4 \}$ is also an edge
of $G_{\cofacet}(P)$. 
\end{enumerate}
In fact, Property~(\ref{item:prop-1}) follows from Property~(\ref{item:prop-2}) 
and the non-existence of triangles, as any larger odd cycle (together with
Property~(\ref{item:prop-2})) implies existence of a triangle. 

We now show that $G_{\cofacet}(P)$ does not have triangles.
Suppose there is a triangle with vertices $v_1, v_2, v_3$ and edges
corresponding to facets $F_1, F_2, F_3$ with $V(F_1) = V(P) \setminus \{v_2,
v_3\}$, $V(F_2) = V(P) \setminus \{ v_1, v_3 \}$, and $V(F_3) = V(P) \setminus
\{ v_1, v_2 \}$.
Let $F'$ be the face $F_1 \cap F_2 = F_1 \cap F_3 = F_2 \cap F_3$. Then
clearly $F_1 = F' * v_1$, $F_2 = F' * v_2$, and $F_3 = F' * v_3$. Thus
$\dim F' = d-2$ and $P/F'$ is a $1$-polytope on $3$-vertices, a contradiction.

Finally, we show that a path $v_1v_2v_3v_4$ of length $3$ implies
the existence of the edge $\{ v_1, v_4 \}$. Let the edges of the path
$v_1v_2v_3v_4$ correspond to facets $F_1, F_2$, and $F_3$ with
$V(F_1) = V(P) \setminus \{ v_1, v_2\}$, $V(F_2) = V(P) \setminus \{ v_2, v_3
\}$, and $V(F_3) = V(P) \setminus \{ v_3, v_4 \}$. 
Let $F' = F_1 \cap F_2 \cap F_3$. Then clearly $F'$ has dimension $d-3$.
Since $F_1$, $F_2$ and $F_3$ are of dimension $d-1$ and each of them contains
exactly two more vertices than $F'$, we conclude that
$F' * v_1$, $F' * v_2$, $F' * v_3$, and $F' * v_4$ are all faces of $P$.
Thus, $P/F'$ is a $2$-polytope on $4$ vertices, which implies that
$F_4 \coloneqq (F' * v_2) * v_3$ is also a facet of $P$ with 
$V(F_4) = V(P) \setminus \{ v_1, v_4\}$.
\end{proof}

\begin{theorem}
\label{thm:polytope-char}
Let $P$ be a $d$-polytope on $d+\gamma+1$ vertices. Then the following are equivalent:
\begin{enumerate}[\rm (i)]
\item Every facet $F$ of $P$ satisfies $|V(P) \setminus V(F)| \leq 2$.
\item $P$ is of the form $\Pnjkm$.
\item $P$ does not have a simplex face of dimension $d-\gamma+1$.
\end{enumerate}
\end{theorem}

\begin{proof}
If $|V(P) \setminus V(F)| \leq 2$ for every facet $F$ of $P$, then
by Lemma~\ref{lem:comb-char} $P$ is of the form $\Pnjkm$.

Now, suppose $P$ is an iterated pyramid over a join of sums of simplices. 
Let $S$ be a simplex face of $P$ of maximal dimension. Then $S$ is the join of
$\Delta_{n-1}$ with facets from each factor $\Delta_{j_i} \oplus \Delta_{k_i}$.
A facet of this sum in turn is obtained by leaving out a vertex from each of
the two simplices. Hence, $S$ has
\[
n + j_1 + k_1 + \ldots + j_m + k_m \; = \; d - m + 1 \; = \; d - \gamma + 1
\]
vertices and therefore dimension $d-\gamma$.

Finally, if $P$ does not have a simplex face of dimension $d-\gamma+1$, then
$|V(P) \setminus V(F)| \leq 2$ for every facet $F$. Otherwise, 
suppose there is a facet $F$ with $|V(P) \setminus V(F)| \geq 3$. 
$\gamma(F) \leq \gamma - 2$, and by Lemma~\ref{lem:simplex} the facet $F$ has a simplex face
of dimension
\[
(d-1) - \gamma(F) = d - (\gamma(F) +1) \geq d - \gamma +1.
\]
\end{proof}

Theorem~\ref{thm:polytope-char} contains the classification of polytopes
on $d+2$ vertices, compare~\cite[pp.~97--101]{Gruenbaum1967}: No $d$-polytope on $d+2$ vertices contains a
simplex $d$-face. Thus, all polytopes on $d+2$ vertices are of
type $\Pnjkm$ with $m=\gamma=1$.

\begin{lemma}
\label{lem:K-subgraph}
Let $P$ be a $d$-polytope on $d+\gamma+1$ vertices. Suppose that
the graph $G(P)$ does not have a $K_{d-\gamma+2}$-subgraph.
Then $P$ is of the form
\[ 
\Pnm = \PnmF,
\]
with $n = d - 3\gamma + 1$ and $m=\gamma$.
\end{lemma}

\begin{proof}
Since $P$ does not have a $K_{d-\gamma+2}$-subgraph, $P$ does
not have a simplex face of dimension $d-\gamma+1$. Thus,
by Theorem~\ref{thm:polytope-char},
$P$ is of the form $P(n, j_1, k_1, \ldots, j_m, k_m)$. 

To show that $j_1 = k_1 = \ldots = j_m = k_m = 1$ observe that the graph
\[
G(\Delta_j \oplus \Delta_k) \left\{ \begin{array}{l@{\text{ if }}l}
    \text{is the complete graph } K_{j+k+2} & j,k \geq 2 \\
    \text{contains a } K_{j+k+1} & j \geq 2, k=1 \text{ or } j=1, k \geq 2 \\
    \text{is a $4$-cycle } & j=k=1 .
\end{array} \right.
\]
Furthermore, in a join $P*Q$ every vertex of $P$ defines an edge with
every vertex of $Q$. Suppose now that $j_i \geq 2$ or $k_i \geq 2$ for some $i$.
Then $G(P)$ contains a complete graph on
\[
n + j_1 + k_1 + \ldots + j_i + k_i + 1 + \ldots + j_m + k_m \; = \; d - m + 2 \; = \; d - \gamma + 2
\]
vertices, but this contradicts the hypothesis.
\end{proof}

\begin{theorem}
\label{thm:char}
Let $P$ be a $d$-polytope on $d+\gamma+1$ vertices
with $k(P) = \left\lfloor (d - \gamma + 1)/2 \right\rfloor$
and $n=d-3\gamma+1$, $m=\gamma$.

If $d-\gamma$ is even, then
\[
P = \Pnm = \PnmF.
\]

If $d - \gamma$ is odd,
there are three possibilities:
\begin{enumerate}[\rm (i)]
\item  \label{item:char-1}
$P = \Pnm$, or
\item \label{item:char-2}
$P = P(n-1,\underbrace{1 ,\ldots, 1}_\text{\parbox{1mm}{\setlength{\unitlength}{1ex}\begin{picture}(0,2)\put(-7,0){$2m-1$ times}\end{picture}}}, 2)$, or
\item \label{item:char-3}
$P$ has a facet $F$ with
\[
F = \Delta_{n-2} * \underbrace{\square * \cdots * \square}_{\text{$m-2$ times}}.
\]
In particular, $k(F) = k(P)$.
\end{enumerate}
\end{theorem}

\begin{proof}
Let $d-\gamma$ be even. If $k(P) = \left\lfloor (d - \gamma + 1)/2 \right\rfloor$, 
then $P$ cannot have a $K_{d-\gamma+2}$-subgraph, and by
Lemma~\ref{lem:K-subgraph} we have $P = \Pnm$ with $n = d - 3\gamma + 1$ and $m=\gamma$.

Let $d-\gamma$ be odd. If $P$ does not have a $K_{d-\gamma+2}$ subgraph, 
then again $P = \Pnm$. So suppose that $P$ does have a $K_{d-\gamma+2}$
subgraph, but not a $K_{d-\gamma+3}$-subgraph. Thus $P$ does not
have a $d-\gamma+2$ simplex face. 
If $P$ also does not have a $d-\gamma+1$ simplex face, then 
\[
P = P(n-1,\underbrace{1 ,\ldots, 1}_{2m-1\text{ times}}, 2).
\]

Consider now the case that $P$ does have a $d-\gamma+1$ simplex face but not a
$d-\gamma+2$ simplex face.
Then for every facet $F$ of $P$ we have $|V(P) \setminus V(F)| \leq 3$. 
By Theorem~\ref{thm:polytope-char} there has to be a facet $F$ with
$|V(P) \setminus V(F)| = 3$, and this facet must satisfy $|V(F) \setminus V(F')| \leq 2$
for every facet $F'$ of $F$. Otherwise there is a ridge $F'$ of $P$ with
a simplex face of dimension at least $d-\gamma+2$.
Then $F$ is of the form $\Pnmp{n'}{m'}$ with $n'=d-3\gamma$ and $m'=\gamma-2$ 
since $d-1 - \gamma(F) = d-1 - (\gamma-2) = d-\gamma+1$ is even and $k(F)=k(P)$.
\end{proof}
The last theorem has interesting consequences.
It implies that for $\gamma > (d+2)/5$ polytopes meeting the
lower bound of $\lfloor (d-\gamma+1)/2 \rfloor$ do not exist. 
Polytopes that appear in the theorem are all at least
$\lfloor (d+\gamma+1)/3 \rfloor$-linked if $\gamma > (d+2)/5$.
But in that case this value is strictly larger than the lower bound.

Furthermore, for $\gamma \leq (d+2)/5 $  and $d-\gamma$ even, 
polytopes meeting the lower bound  
are unique. Thus, they are characterized by
Theorem~\ref{thm:char}. 

However, if $d-\gamma$ is odd, such polytopes are not characterized
by the three possibilities given.
While the polytopes in Possibility~(\ref{item:char-1}) and~(\ref{item:char-2})
are $\lfloor (d-\gamma+1)/2 \rfloor$-linked, polytopes
as in Possibility~(\ref{item:char-3}) can be higher linked.
We find different examples of type~(\ref{item:char-3}) by
replacing certain factors of the join in $\Pnm$. 

If we replace the $5$-dimensional polytope
$Q := \square * \square$ by the two-fold pyramid over the $3$-dimensional
crosspolytope, that is, 
\[
\square * \square \quad \rightsquigarrow \quad \Delta_1 * C_3^\Delta, 
\]
we obtain a polytope $P$ with $k(P) = k(\Pnm)$:
The polytope $\Delta_1 * C_3^\Delta$ is $5$-dimensional and has the
same number of vertices as $Q$. 
In the complement of the graph $G(\Pnm)$ four isolated edges are replaced
by three isolated edges and two isolated vertices.
This change does not increase linkedness, as the
condition that $d-\gamma$ is odd in terms of $n$ and $m$
translates to the condition that $n$ is even. Hence, we have $k(P) = k(\Pnm)$.

Similar observations show that if we replace $Q$ by the
two-fold pyramid over a triangular prism we also obtain a 
polytope $P$ with $k(P) = k(\Pnm)$.

However, it is possible for a polytope to have a facet as in
Possibility~(\ref{item:char-3}) and nevertheless to be higher linked
than $\Pnm$. We obtain such a polytope $P$ for instance if we replace
the factor $Q$ in $\Pnm$ by a two-fold pyramid over a twice stacked $3$-simplex.

\section{Conclusions and open problems}

Theorem~\ref{thm:general-lower} and Corollary~\ref{cor:gallivan-bound} 
imply the values for $k(d)$ as displayed in Table~\ref{tab:tabular-kd}.

\begin{table}[hb]
  \centering
  \begin{minipage}[h]{.22\textwidth}
    \centering
    \begin{tabular}{rr}
      $d$ & ~~$k(d)$ \\
      \hline
      $1$ & 1 \\
      $2$ & 1 \\
      $3$ & 1 \\
      $4$ & 2 \\
      $5$ & 2
    \end{tabular}
  \end{minipage}
  \begin{minipage}[h]{.22\textwidth}
    \centering
    \begin{tabular}{rr}
      $d$ & ~~$k(d)$ \\
      \hline
      6 & 2,3 \\
      7 & 3 \\
      8 & 3 \\
      9 & 3,4 \\
      10 & 4 
    \end{tabular}
  \end{minipage}
  \begin{minipage}[h]{.22\textwidth}
    \centering
    \begin{tabular}{rr}
      $d$ & ~~$k(d)$ \\
      \hline
      11 & 4,5 \\
      12 & 4,5 \\
      13 & 5 \\
      14 & 5,6 \\
      15 & 5,6,7
    \end{tabular}
  \end{minipage}
  \caption{Possible values of $k(d)$.}
  \label{tab:tabular-kd}
\end{table}

 In particular, we
 get exact values in dimensions $7$, $10$, and $13$. The value $k(8) =3$
 follows from Larman and Mani's old lower bound~\cite{LarmanMani1970} and
 Gallivan's upper bound~\cite{Gallivan1985}.

 The first open value is $k(6)$ and it seems to be a difficult problem to
 determine it. Our analysis  of polytopes with few vertices (Theorem~\ref{thm:char}) shows
 that $k(6,0) = k(6,1) = k(6,2)  = 3$. We have also verified enumeratively
 that $k(6,3) = 3$; beyond that we do not know anything.

 \begin{prob}
 Determine $k(6)$: Either show that all $6$-polytopes are $3$-linked, 
 or give an example of a $6$-polytope $P$ with $k(P) = 2$.
 \end{prob}

 One can construct polytopes with $f_0 = 3 \lfloor d/2 \rfloor - 1$
 vertices that are not $\lfloor d/2 \rfloor$-linked, which is the
 bound in the original question by Larman and Mani.
 If $d$ is even let 
 \[
 P \; := \; \Delta_2 * \square * \square * \underbrace{C^\Delta_3 * \cdots *
   C^\Delta_3}_{\text{$m$ times}}.
 \]
 Then $d = 4m + 8$, $f_0 = 6m + 11$ and $k(P) = 2m + 3$.

 For $d$ odd let 
 \[
 P \; := \; \Delta_4 * \square * \square * \square * \underbrace{C^\Delta_3 * \cdots *
   C^\Delta_3}_{\text{$m$ times}}.
 \]
 Then $d = 4m+13$, $f_0 = 6m+17$ and $k(P) = 2m + 5$.

 \begin{prob}
   \label{prob:many-vertices}
 Are all $d$-polytopes on at least $3 \lfloor d/2 \rfloor$ vertices $\lfloor d/2
 \rfloor$-linked?
 Weaker: Is there some $N(d)$, such that every $d$-polytope on
 at least $N(d)$ vertices is $\lfloor d/2 \rfloor$-linked?
 \end{prob}

\begin{figure}[ht]
  \centering
  \includegraphics[width=.2\textwidth]{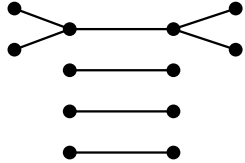}
  \caption{Is this a subgraph of the complement graph of some $8$-polytope
    on $12$ vertices?}
  \label{fig:8-graph}
\end{figure}

 Only one obstruction for $d$-polytopes to not be $\lfloor d/2 \rfloor$-linked
 is known, the obstruction exploited by Gallivan: The polytopes have
 many missing edges and not enough vertices to route all paths around
 the missing edges.
 If a polytope has $3 \lfloor d/2 \rfloor$ or more vertices, there has to be
 a different obstruction if it is not $\lfloor d/2 \rfloor$-linked.
 Regarding Problem~\ref{prob:many-vertices}, it would be interesting to know if the graph
 in Figure~\ref{fig:8-graph}, which is not $4$-linked, is a subgraph of the
 complement graph of an $8$-polytope on $12$ vertices. The complement of this
 graph is $8$-connected, and at every vertex it has a subdivision of $K_9$ rooted at that vertex.

\vfill

\bibliographystyle{siam}

\end{document}